\documentstyle[amscd,amssymb,verbatim,12pt]{amsart}
\newcommand{\Cobar}{\operatorname{Cobar}}
\renewcommand{\Bar}{\operatorname{Bar}}
\newcommand{\Fun}{\operatorname{Fun}}
\newcommand{\gr}{\operatorname{gr}}
\newcommand{\intdeg}{\operatorname{intdeg}}

\newcommand{\OO}{{\cal O}}

\newcommand{\limproj}{\operatorname{proj.lim}}

\newcommand{\DD}{{\cal D}}

\newcommand{\lan}{\langle}
\newcommand{\ran}{\rangle}

\newcommand{\CC}{{\cal C}}

\newcommand{\Spec}{\operatorname{Spec}}

\newcommand{\ga}{\gamma}
\newcommand{\de}{\delta}
\newcommand{\eps}{\epsilon}

\renewcommand{\ker}{\operatorname{ker}}

\numberwithin{equation}{subsection}

\newtheorem{thm}{Theorem}[section]
\newtheorem{prop}[thm]{Proposition}
\newtheorem{lem}[thm]{Lemma}
\newtheorem{cor}[thm]{Corollary}
\newenvironment{rem}{\vspace{3mm}\noindent
{\bf Remark.}}{\vspace{3mm}}
\newenvironment{rems}{\vspace{3mm}
\noindent {\bf Remarks.}}{\vspace{3mm}}

\newcommand{\Pf}{\noindent {\it Proof}}
\newcommand{\id}{\operatorname{id}}

\newcommand{\ra}{\rightarrow}

\newcommand{\FF}{{\cal F}}

\newcommand{\Hom}{\operatorname{Hom}}
\newcommand{\Ext}{\operatorname{Ext}}

\renewcommand{\a}{\alpha}
\renewcommand{\b}{\beta}

\newcommand{\la}{\lambda}

\newcommand{\C}{{\Bbb C}}
\newcommand{\A}{{\Bbb A}}

\newcommand{\Z}{{\Bbb Z}}

\newcommand{\wt}{\widetilde}

\newcommand{\sub}{\subset}
\newcommand{\ed}{\qed\vspace{3mm}}

\def\cech{\v Cech}

\title{Extensions of homogeneous coordinate rings to $A_{\infty}$-algebras}
\author{A. Polishchuk}
\thanks{This work was partially supported by NSF grant DMS-0070967}

\begin{document}
\begin{abstract}
We study $A_{\infty}$-structures extending the natural algebra structure on
the cohomology of $\oplus_{n\in\Z} L^n$,
where $L$ is a very ample line bundle on a projective $d$-dimensional
variety $X$ such that $H^i(X,L^n)=0$ for $0<i<d$ and all $n\in\Z$. 
We prove that there exists a unique such {\it nontrivial}
$A_{\infty}$-structure up to a strict $A_{\infty}$-isomorphism (i.e., an
$A_{\infty}$-isomorphism with the identity as the first structure map)
and rescaling.
In the case when $X$ is a curve we also compute the group of 
strict $A_{\infty}$-automorphisms of this $A_{\infty}$-structure.
\end{abstract}
\maketitle

\bigskip

\section{Introduction}

Let $X$ be a projective variety over a field $k$, 
$L$ be a very ample line bundle on $X$. Recall that
the graded $k$-algebra
$$R_L=\oplus_{n\ge 0}H^0(X,L^n)$$
is called the {\it homogeneous coordinate ring}
corresponding to $L$. More generally, one can consider 
the bigraded $k$-algebra
$$A_L=\oplus_{p,q\in\Z} H^q(X,L^p).$$
We call $A_L$ the {\it extended homogeneous coordinate ring}
corresponding to $L$. 

Since $A_L$ can be represented naturally as the cohomology algebra of some
dg-algebra (say, using injective resolutions or
{\cech} cohomology with respect to an affine
covering), it is equipped with a family of higher operations called
Massey products. A better way of recording this additional structure 
uses the notion of $A_{\infty}$-algebra due to Stasheff. Namely, by 
the theorem of Kadeishvili the product on $A_L$ extends to
a canonical (up to $A_{\infty}$-isomorphism) $A_{\infty}$-algebra structure
with $m_1=0$ (see \cite{Keller} 3.3 and references therein). 
More precisely, this structure is unique up to a {\it strict
$A_{\infty}$-isomorphism}, i.e., an $A_{\infty}$-isomorphism with 
the identity map as the first structure map (see section \ref{strsec} for
details).
Note that the axioms of $A_{\infty}$-algebra use the
{\it cohomological} grading on $A_L$ (where $H^q(X,L^p)$ has 
cohomological degree $q$),
and all the operations $(m_n)$ have degree zero with
respect to the {\it internal} grading
(where $H^q(X,L^p)$ has {\it internal degree} $p$).
The natural question is whether it is possible
to characterize intrinsically this canonical class of
$A_{\infty}$-structures on $A_L$. This question is partly motivated by
the homological mirror symmetry. Namely, in the case when $X$ is a
Calabi-Yau manifold, the $A_{\infty}$-structure on $A_L$
is supposed to be $A_{\infty}$-equivalent to an
appropriate $A_{\infty}$-algebra
arising on a mirror dual symplectic side. An intrinsic characterization of
the $A_{\infty}$-isomorphism class of our 
$A_{\infty}$-structure could be helpful
in reducing the problem of constructing such an $A_{\infty}$-equivalence
to constructing an isomorphism of the usual associative algebras.
More generally, it is conceivable that the algebra $A_L$ can appear
as cohomology algebra of some other dg-algebras (for example, if
there is an equivalence of the derived category of coherent
sheaves on $X$ with some other such category), so one might be
interested in comparing corresponding $A_{\infty}$-structures
on $A_L$.

Thus, we want to study all $A_{\infty}$-structures $(m_n)$ on $A_L$
(with respect to the cohomological grading),
such that $m_1=0$, $m_2$ is the standard double
product and all $m_n$ have degree $0$ with respect to the
internal grading. Let us call such an $A_{\infty}$-structure
on $A_L$ {\it admissible}. As we have already mentioned before,
there is a {\it canonical strict $A_{\infty}$-isomorphism class}
of such structures coming from the realization of $A_L$ as cohomology
of the dg-algebra $\oplus_n \CC^{\bullet}(L^n)$
where $\CC^{\bullet}(?)$ denotes the {\cech} complex with respect to some
open affine covering of $X$. By definition,
an $A_{\infty}$-structure belongs to the canonical class if there exists
an $A_{\infty}$-morphism from $A_L$ 
equipped with this $A_{\infty}$-structure to the above dg-algebra
inducing identity on the cohomology. The simplest picture one could imagine
would be that all admissible $A_{\infty}$-structures are 
strictly $A_{\infty}$-isomorphic, i.e., that $A_L$ is intrinsically formal. 
It turns out that
this is not the case. However, our main theorem
below shows that if the cohomology of the structure sheaf on $X$ is
concentrated in degrees $0$ and $\dim X$ then
for sufficiently ample $L$ the situation is
not too much worse.

We will recall the notion of a homotopy between $A_{\infty}$-morphisms
in section \ref{strsec} below. \footnote{All $A_{\infty}$-morphisms
and homotopies between them are assumed to respect the
internal grading on $A_L$.}
Let us say that an $A_{\infty}$-structure is {\it nontrivial} if it is not
$A_{\infty}$-isomorphic to an 
$A_{\infty}$-structure with $m_n=0$ for $n\neq 2$.
By rescaling of an $A_{\infty}$-structure we mean the
change of the products $(m_n)$ to $(\la^{n-2}m_n)$ for some
constant $\la\in k^*$.
Our main result gives a classification of admissible $A_{\infty}$-structures
on $A_L$ up to a strict $A_{\infty}$-isomorphism and 
rescaling (under certain assumptions).

\begin{thm}\label{mainthm}
Let $L$ be a very ample line bundle on a $d$-dimensional
projective variety $X$ such that
$H^q(X,L^p)=0$ for $q\neq 0,d$ and all $p\in\Z$. 
Then

\noindent
(i) up to a strict $A_{\infty}$-isomorphism
and rescaling there exists a unique non-trivial
admissible $A_{\infty}$-structure on $A_L$;
moreover, $A_{\infty}$-structures on $A_L$ 
from the canonical strict $A_{\infty}$-isomorphism class are nontrivial;

\noindent
(ii) for every pair of strict $A_{\infty}$-isomorphisms 
$f,f':(m_i)\ra(m'_i)$ between
admissible $A_{\infty}$-structures on $A_L$ there exists a homotopy
from $f$ to $f'$.
\end{thm}

\begin{rems} 1. One can unify strict $A_{\infty}$-isomorphisms
with rescalings by considering $A_{\infty}$-isomorphisms $(f_n)$
with the morphism $f_1$ of the form $f_1(a)=\la^{\deg(a)}$
for some $\la\in k^*$ (where $a$ is a homogeneous element of $A_L$).
In particular, part (i) of the theorem implies that all non-trivial
admissible $A_{\infty}$-structures on $A_L$ are $A_{\infty}$-isomorphic.

\noindent
2. As we will show in section \ref{strsec}, part (ii) of the theorem
is equivalent to its particular case when $f'=f$. In this case the statement
is that every strict $A_{\infty}$-automorphism of $f$ is homotopic to
the identity.

\noindent
3. If one wants to see more explicitly how a canonical $A_{\infty}$-structure on $A_L$
looks like, one has to choose one of the natural dg-algebras with cohomology $A_L$ (an
obvious algebraic choice is the {\cech} complex; in the case $k=\C$ one can also use the
Dolbeault complex), choose a projector $\pi$ from the dg-algebra to some 
space of representatives for the cohomology such that $\pi=1-dQ-Qd$ for some
operator $Q$, and then apply formulas of \cite{KS} for the operations $m_n$
(they are given by certain sums over trees).  
\end{rems}

The above theorem is applicable to every
line bundle of sufficiently large degree on a curve.
In higher dimensions it can be used for every sufficiently
ample line bundle on a $d$-dimensional projective variety $X$ such that 
there exists a dualizing sheaf on $X$ and
$H^i(X,\OO_X)=0$ for $0<i<d$. For example, this condition
is satisfied for complete intersections in projective spaces.
At present we do not know how
to extend this theorem to the case when $\OO_X$ has some
nontrivial middle cohomology. Note that for a smooth
projective variety $X$ over $\C$ the natural (up to a strict 
$A_{\infty}$-isomorphism)
$A_{\infty}$-structure on $H^*(X,\OO_X)$ is trivial as follows
from the formality theorem of \cite{DGMS}. This suggests
that for sufficiently ample line bundle $L$ 
one could try to characterize the canonical $A_{\infty}$-structure
on $A_L$ (up to a strict $A_{\infty}$-isomorphism and rescaling) 
as an admissible
$A_{\infty}$-structure whose restriction to $H^*(X,\OO_X)$ 
is trivial.

In the case when $X$ is a curve we can also compute the group of
strict $A_{\infty}$-automorphisms of an $A_{\infty}$-structure on $A_L$.
As we will explain below \ref{strsec}, strict $A_{\infty}$-isomorphisms
on $A_L$ form a group $HG$, which is a subgroup of 
automorphisms of the free coalgebra
$\Bar(A_L)$ (preserving both gradings).
The dual to the degree zero component 
of $\Bar(A_L)$ (with respect to both gradings)
can be identified with the completed
tensor algebra $\hat{T}(H^1(X,\OO_X)^*)=\prod_{n\ge 0}T^n(H^1(X,\OO_X)^*)$. 
Therefore, we obtain a natural homomorphism from $HG$ to 
to the group $G$ of continuous automorphisms of 
$\hat{T}(H^1(X,\OO_X)^*)$. 

\begin{thm}\label{groupthm}
Let $L$ be a very ample line bundle on a projective curve $X$ 
such that $H^1(X,L)=0$. Let also $HG(m)\sub HG$ be the group of 
strict $A_{\infty}$-automorphisms of an admissible 
$A_{\infty}$-structure $m$ on $A_L$.
Then the above homomorphism $HG\ra G$ restricts to an isomorphism of 
$HG(m)$ with the subgroup $G_0\sub G$ consisting of inner automorphisms 
of $\hat{T}(H^1(X,\OO_X)^*)$ by elements in
$1+\prod_{n>0}T^n(H^1(X,\OO_X)^*)$. 
\end{thm}

Assume that $X$ is a smooth projective curve. Then there is a canonical
noncommutative thickening $\wt{J}$ of the Jacobian $J$ of $X$ (see
\cite{Kap}). As was shown in \cite{P-b}, 
a choice of an $A_{\infty}$-structure in the canonical strict 
$A_{\infty}$-isomorphism
class gives rise to a formal system of coordinates on $\wt{J}$ at zero.
More precisely, by this we mean an isomorphism of the formal
completion of the local ring of $\wt{J}$ at zero with $\hat{T}(H^1(X,\OO_X)^*)$
inducing the identity map on the tangent spaces.
Formal coordinates associated with two 
strictly isomorphic $A_{\infty}$-structures 
are related by the coordinate change given
by the image of the corresponding $A_{\infty}$-isomorphism 
under the homomorphism $HG\ra G$. 
Now Theorem \ref{groupthm}
implies that two $A_{\infty}$-structures in the canonical 
class that induce the same formal coordinate on $\wt{J}$ can be connected
by a unique strict $A_{\infty}$-isomorphism. 
Indeed, two such isomorphisms would differ by a 
strict $A_{\infty}$-automorphism
inducing the trivial automorphism of $\hat{T}(H^1(X,\OO_X)^*)$,
but such an $A_{\infty}$-automorphism is trivial by Theorem \ref{groupthm}.

\vspace{2mm}

\noindent
{\it Convention}. Throughout the paper we work over a fixed ground field $k$.
The symbol $\otimes$ without additional subscripts always denotes the tensor
product over $k$.

\vspace{2mm}

\noindent
{\it Acknowledgment}. I'd like to thank the referee for helpful 
remarks and suggestions.

\section{Preliminaries}

\subsection{Strict $A_{\infty}$-isomorphisms and homotopies}\label{strsec}

We refer to \cite{Keller} for an introduction to
$A_{\infty}$-structures. We restrict ourselves to several remarks
about $A_{\infty}$-isomorphisms and homotopies between them.

A {\it strict $A_{\infty}$-isomorphism} between two $A_{\infty}$-structures
$(m)$ and $(m')$ on the same graded space $A$ is an
$A_{\infty}$-morphism $f=(f_n)$ from $(A,m)$ to $(A,m')$ such that
$f_1=\id$. The equations connecting $f$, $m$ and $m'$ can be interpreted
as follows.
Recall that $m$ and $m'$ correspond to coderivations $d_m$ and $d_{m'}$
of the bar-construction $\Bar(A)=\oplus_{n\ge 1} T^n(A[1])$ such that
$d_m^2=d_{m'}^2=0$. Now every collection $f=(f_n)_{n\ge 1}$, where 
$f_n:A^{\otimes n}\ra A$ has degree $1-n$, $f_1=\id$, defines
a coalgebra automorphism $\a_f:\Bar(A)\ra \Bar(A)$, 
with the component $\Bar(A)\ra A[1]$ given by $(f_n)$.
The condition that $f$ is an $A_{\infty}$-morphism is equivalent
to the equation $\a_f\circ d_m= d_{m'}\circ \a_f$. 
In other words, 
strict $A_{\infty}$-isomorphisms between $A_{\infty}$-structures precisely
correspond to the action of the group of automorphisms of $\Bar(A)$
as a coalgebra on the space of coderivations $d$ such that
$d^2=0$. More precisely, we consider only automorphisms of $\Bar(A)$
of degree $0$ inducing the identity map $A\ra A$.
Let us denote by
$HG=HG(A)$ the group of such automorphisms which we will also call
{\it the group of strict $A_{\infty}$-isomorphisms} on $A$.
We will denote by $m \to g*m$, where $g\in HG$, the natural action
of this group on the set of all $A_{\infty}$-structures on $A$.

One can define a decreasing filtration
$(HG_n)$ of $HG$ by normal subgroups by setting 
$$HG_n=\{ f=(f_i)\ |\ f_i=0, 1<i\le n\}.$$ 
Note that $f\in HG_n$ if and only if the restriction of
$\a_f$ to the sub-coalgebra
$\Bar(A)_{\le n}=\oplus_{i\le n}(A[1])^{\otimes i}$
is the identity homomorphism. Furthermore, it is also clear
that $HG\simeq \projlim_n HG/HG_n$. In particular, an infinite
product of strict $A_{\infty}$-isomorphisms 
$\ldots *f(3)*f(2)*f(1)$ is well-defined as long
as $f(n)\in HG_{i(n)}$, where $i(n)\to\infty$ as $n\to\infty$.

The notion of a homotopy between $A_{\infty}$-morphisms is best understood
in a more general context of $A_{\infty}$-categories. Namely,
for every pair of $A_{\infty}$-categories $\CC$, $\DD$ one can define
the $A_{\infty}$-category $\Fun(\CC,\DD)$ having 
$A_{\infty}$-functors $F:\CC\ra\DD$ as objects (see \cite{Lefevre},
\cite{Lyub}).
In particular, there is a natural notion of closed morphisms between
two $A_{\infty}$-functors $F,F':\CC\ra\DD$. 
Specializing to the case when $\CC$ and $\DD$ are 
$A_{\infty}$-categories with one
object corresponding to $A_{\infty}$-algebras $A$ and $B$ 
we obtain a notion of a closed morphism between
a pair of $A_{\infty}$-morphisms $f,f':A\ra B$. Following
\cite{Keller} we call such a closed morphism a 
{\it homotopy between $A_{\infty}$-morphisms $f$ and $f'$}.
More explicitly, a homotopy $h$ is given by a 
collection of maps $h_n:A^{\otimes n}\ra B$
of degree $-n$, where $n\ge 1$, satisfying some equations.
These equations are written as follows: there exists a unique linear map
$H:\Bar(A)\ra\Bar(B)$ of degree $-1$ with the component
$\Bar(A)\ra B$ given by $(h_n)$, such that
\begin{equation}\label{hommor1}
\Delta\circ H=(\a_f\otimes H+H\otimes \a_{f'})\circ\Delta,
\end{equation}
where $\a_f,\a_{f'}:\Bar(A)\ra\Bar(B)$ are coalgebra
homomorphisms corresponding to $f$ and $f'$, $\Delta$ denotes
the comultiplication.
Then the equation connecting $h$, $f$ and $f'$ is
\begin{equation}\label{hommor2}
\a_f-\a_{f'}=d_A\circ H+H\circ d_B,
\end{equation}
where $d_A$ (resp., $d_B$) is the coderivation of $\Bar(A)$ (resp., $\Bar(B)$)
corresponding to the $A_{\infty}$-structure on $A$ (resp., $B$).
It is not difficult to check that for a given $A_{\infty}$-morphism
$f$ from $A$ to $B$ the equations (\ref{hommor1}) and (\ref{hommor2})
imply that $\a_{f'}$ is a homomorphism of dg-coalgebras, so it defines
an $A_{\infty}$-morphism $f'$ from $A$ to $B$. Moreover, similarly
to the case of strict $A_{\infty}$-isomorphisms we have the following result.

\begin{lem}\label{closmorlem} 
Let $A$ and $B$ be $A_{\infty}$-algebras and $f=(f_n)$
be an $A_{\infty}$-morphism from $A$ to $B$. For every
collection $(h_n)_{n\ge 1}$, where $h_n:A^{\otimes n}\ra B$ has
degree $-n$, there exists a unique $A_{\infty}$-morphism $f'$ from
$A$ to $B$ such that $h$ is a homotopy from $f$ to $f'$.
\end{lem}

\Pf . It is easy to see that equation (\ref{hommor1}) is
equivalent to the following formula
\begin{equation}\label{hommor3}
\begin{array}{l} H[a_1|\ldots|a_n]=\sum_{i_1<\ldots<i_k<m<j_1<\ldots<j_l=n}
\pm [f_{i_1}(a_1,\ldots,a_{i_1})|\ldots 
|f_{i_k-i_{k-1}}(a_{i_{k-1}+1},\ldots,a_{i_k})|  \\
h_{m-i_k}(a_{i_k+1},\ldots a_m)|f'_{j_1-m}(a_{m+1},\ldots,a_{j_1})|\ldots|
f'_{j_l-j_{l-1}}(a_{j_{l-1}+1},\ldots,a_{j_l})],
\end{array}
\end{equation}
where $a_1,\ldots,a_n\in A$, $n\ge 1$.
We are going to construct the maps $H|_{\Bar(A)_{\le n}}$ and
$\a_{f'}|_{\Bar(A)_{\le n}}$ recursively, so that at each
step the equations (\ref{hommor2}) and (\ref{hommor3})
are satisfied when restricted to $\Bar(A)_{\le n}$.
Of course, we also want $H$ to have $(h_n)$ as components.
Then such a construction will be unique.
Note that $H|_{A[1]}$ is given by $h_1$ and
$\a_{f'}|_{A[1]}$ is given by $f'_1=f_1-m_1\circ h_1-h_1\circ m_1$.
Now assume that the restrictions of $H$ and $\a_{f'}$
to $\Bar(A)_{\le n-1}$ are already constructed, so that the maps
$f'_i:A^{\otimes i}\ra B$ are defined for $i\le n-1$.
Then the formula (\ref{hommor3}) defines uniquely the extension of
$H$ to $\Bar(A)_{\le n}$ (note that in the RHS of this formula
only $f'_i$ with $i\le n-1$ appear). It remains to apply formula
(\ref{hommor2}) to define $\a_{f'}|_{\Bar(A)_{\le n}}$.
\ed

Let $HG$ be the group of strict $A_{\infty}$-isomorphisms on a given
graded space $A$. In other words, $HG$ is the group of
degree $0$ coalgebra automorphisms of $\Bar(A)$ with the componenent
$A\ra A$ equal to the identity map. This group acts on the
set of all $A_{\infty}$-structures on $A$. The stabilizer subgroup
of some $A_{\infty}$-structure $m$ is the group of strict 
$A_{\infty}$-automorphisms $HG(m)$.
We can consider the set of all strict $A_{\infty}$-automorphisms $f_h\in HG(m)$
such that there exists a homotopy $h$ from the trivial 
$A_{\infty}$-automorphism
$f^{tr}$ to $f_h$ (where $f^{tr}_i=0$ for $i>1$). It is easy to see that
$A_{\infty}$-automorphisms of the form $f_h$ constitute a normal subgroup
in $HG(m)$ that we will denote by $HG(m)^0$. Furthermore, for every
$g\in HG$ we have $HG_{g*m}=g HG(m)^0 g^{-1}$. Also, for
a pair of elements $g_1,g_2\in HG$ such that $m'=g_1*m=g_2*m$,
there exists a homotopy between $g_1$ and $g_2$ (where
$g_i$ are considered as $A_{\infty}$-morphisms from $(A,m)$ to $(A,m')$)
if and only if $g_1^{-1}g_2\in HG(m)^0$.

\subsection{Obstructions}\label{obstrsec}

Below we use Hochschild cohomology $HH(A):=HH(A,A)$
for a graded associative algebra $A$
(see \cite{Lod} for the corresponding sign convention).
When considering $A=A_L$ as a graded algebra we
equip it with the cohomological grading, so in the situation of
Theorem \ref{mainthm} this grading has only $0$-th and $d$-th
non-trivial graded components.

The following lemma is well known and its proof is straightforward.

\begin{lem}\label{obstrlem}
Let $m$ and $m'$ be two admissible
$A_{\infty}$-structures on $A$. Assume that $m_i=m'_i$ for $i<n$,
where $n\ge 3$.

\noindent
(i) Set
$c(a_1,\ldots,a_n)=(m'_n-m_n)(a_1,\ldots,a_n)$.
Then $c$ is a Hochschild $n$-cocycle, i.e.,
\begin{align*}
&\de c(a_1,\ldots,a_{n+1})=\sum_{j=1}^n 
(-1)^j c(a_1,\ldots,a_ja_{j+1},\ldots,a_{n+1})+\\
&(-1)^{n\deg(a_1)}a_1c(a_2,\ldots,a_{n+1})+(-1)^{n+1}
c(a_1,\ldots,a_n)a_{n+1}=0.
\end{align*}

\noindent
(ii) If $m'=f*m$ for a strict $A_{\infty}$-isomorphism $f$ such that $f_i=0$
for $1<i<n-1$, then setting
$b(a_1,\ldots,a_{n-1})=(-1)^{n-1}f_{n-1}(a_1,\ldots,a_{n-1})$
we get
$$c(a_1,\ldots,a_n)=\delta b(a_1,\ldots,a_n),$$
where $c$ is the $n$-cocycle defined in (i). Hence, $c$ is
a Hochschild coboundary.
\end{lem}

Thus, the study of admissible $A_{\infty}$-structures on $A$
is closely related to the study of certain components of
Hochschild cohomology of $A$. More precisely, let us denote
$C^n_{p,q}(A)$ (resp. $HH^n_{p,q}(A)$)
the space of reduced Hochschild $n$-cochains (resp. of $n$-th
Hochschild cohomology classes)
of internal grading $-p$ and of cohomological grading $-q$.          
In other words, $C^n_{p,q}(A)$ consists of cochains
$c:A^{\otimes n}\ra A$ such that
$\intdeg c(a_1,\ldots,a_n)=\intdeg a_1+\ldots+\intdeg a_n-p$,
$\deg c(a_1,\ldots,a_n)=\deg a_1+\ldots+\deg a_n-q$.
Since, all the operations $m_n$ respect the internal
grading and have (cohomological) degree $2-n$, we see that
the cocycle $c$ defined in Lemma \ref{obstrlem}
lives in $C^n_{0,n-2}(A)$.
                                               
There is an analogue of Lemma \ref{obstrlem} for 
strict $A_{\infty}$-isomorphisms.

\begin{lem}\label{obstrlem2}
Let $m$ and $m'$ be admissible $A_{\infty}$-structures on $A$,
$f,f'$ be a pair of strict $A_{\infty}$-isomorphisms 
from $m$ to $m'$. Assume that $f_i=f'_i$
for $i<n$, where $n\ge 2$.

\noindent
(i)  Set
$c(a_1,\ldots,a_n)=(f'_n-f_n)(a_1,\ldots,a_n)$.
Then $c$ is a Hochschild $n$-cocycle in $C^n_{0,n-1}(A)$.

\noindent
(ii) If $\phi:f\ra f'$ is a homotopy such that
$\phi_i=0$ for $i<n-1$, then for
$b(a_1,\ldots,a_{n-1})=\pm \phi_{n-1}(a_1,\ldots,a_{n-1})$
one has $c=\delta b$.
\end{lem}

\section{Calculations}

\subsection{Hochschild cohomology}

In this subsection we calculate the components of
the Hochschild cohomology of $A=A_L$ that are relevant
for the proof of Theorem \ref{mainthm}.

Let us set $R=R_L$ and let $R_+=\oplus_{n\ge 1}R_n$
be the augmentation ideal in $R$, so that $R/R_+=k$.
Recall that the bar-construction provides a
free resolution of $k$ as $R$-module of the form
\begin{equation}\label{bar}
\ldots\ra R_+\otimes R_+\otimes R\ra R_+\otimes R\ra R\ra k.
\end{equation}
For graded $R$-bimodules $M_1,\ldots, M_n$ we consider the bar-complex
$$B^{\bullet}(M_1,\ldots,M_n)=M_1\otimes T(R_+)\otimes M_2\otimes
\ldots T(R_+)\otimes M_n,$$
where $T(R_+)$ is the tensor algebra of $R_+$. 
This is just the tensor product over $T(R_+)$ of the bar-complexes
of $M_1,\ldots, M_n$ (where $M_1$ is considered as a right $R$-module,
$M_2,\ldots,M_{n-1}$ as $R$-bimodules, and $M_n$ as a left $R$-module). 
The grading in this complex
is induced by the {\it cohomological grading} of the tensor algebra
$T(R_+)$ defined by $\deg T^i(R_+)=-i$,
so that $B^{\bullet}(M_1,\ldots,M_n)$ is concentrated in nonnegative
degrees and the differential has degree $1$.
For example, $B^{\bullet}(k,R)$ is the bar-resolution (\ref{bar}) of $k$.

\begin{prop}\label{cohprop}
Under the assumptions of Theorem \ref{mainthm}
let us consider the graded $R$-module $M=\oplus_{i\in\Z} H^d(X,L^i)$. 
Let $M_1,\ldots,M_n$ be graded $R$-bimodules such that
each of them is isomorphic to $M$ as a (graded) right $R$-module
and as a left $R$-module.

\noindent
(i) The complex $B^{\bullet}(k,M)=T(R_+)\otimes M$ 
has one-dimensional cohomology,
which is concentrated in degree $-d-1$ and internal degree $0$.

\noindent
(ii) $H^i(B^{\bullet}(M_1,M_2))=0$ for
$i\neq -d-1$ and $H^{-d-1}(B^{\bullet}(M_1,M_2))$ is isomorphic to $M$
as a (graded) right $R$-module and as a left $R$-module. 

\noindent
(iii) $H^i(B^{\bullet}(M_1,\ldots,M_n))=0$ for $i>-(n-1)(d+1)$.

\noindent
(iv) $H^i(B^{\bullet}(k,M_1,\ldots,M_n,k))=0$ for $i>-n(d+1)$.
In addition, the space \break
$H^{-d-1}(B^{\bullet}(k,M_1,k))$ is one-dimensional
and has internal degree $0$.
\end{prop}

\Pf . (i) Localizing the exact sequence (\ref{bar}) on $X$ and tensoring
with $L^m$, where $m\in\Z$, we obtain the
following exact sequence of vector bundles on $X$:
\begin{equation}\label{barsh}
\ldots \oplus_{n_1,n_2>0}R_{n_1}\otimes R_{n_2} \otimes L^{m-n_1-n_2}\ra
\oplus_{n>0} R_n\otimes L^{m-n}\ra L^m\ra 0.
\end{equation}
Each term in this sequence is a direct sum of a number of copies
of line bundles $L^n$: 
for a finite-dimensional vector space $V$ 
we denote by $V\otimes L^n$ the direct sum of $\dim V$ copies of $L^n$. 
Now let us consider the spectral sequence with $E_1$-term given by
the cohomology of all sheaves in this complex and abutting
to zero (this sequence converges since we can compute
cohomology using {\cech} resolutions with respect to a finite open affine
covering of $X$).
The $E_1$-term consists of two rows: one obtained by applying
the functor $H^0(X,\cdot)$ to (\ref{barsh}), another obtained by
applying $H^d(X,\cdot)$. The row of $H^0$'s has form
$$\ldots \oplus_{n_1,n_2>0}R_{n_1}\otimes R_{n_2}\otimes R_{m-n_1-n_2}\ra
\oplus_{n>0} R_n\otimes R_{m-n}\ra R_m\ra 0$$
which is just the $m$-th homogeneous component of the bar-resolution.
Hence, this complex is exact for $m\neq 0$. Since the sequence abuts to
zero the row of $H^d$'s should also be exact for $m\neq 0$. For
$m=0$ the row of $H^0$'s reduces to the single term $H^0(X,\OO_X)=k$,
hence, the row of $H^d$'s in this case has one-dimensional
cohomology at $-(d+1)$th term and is exact elsewhere.

\noindent
(ii) Consider the filtration on $B^{\bullet}(M_1,M_2)$ associated with the
$\Z$-grading on $M_2$. By part (i) the corresponding spectral sequence
has the term $E_1\simeq H^{-d-1}(M_1\otimes T(R_+))\otimes M_2\simeq M_2$.
Hence, it degenerates in this term and
$$H^*(K^{\bullet})\simeq H^{-d-1}(K^{\bullet})\simeq M$$
as a right $R$-module. Similarly, the spectral sequence
associated with the filtration
on $K^{\bullet}$ induced by the $\Z$-grading on $M_2$ 
gives an isomorphism of left $R$-modules $H^{-d-1}(K^{\bullet})\simeq M$.

\noindent
(iii) For $n=2$ this follows from (ii). Now let $n>2$ and assume that the
assertion holds for $n'<n$. We can consider $B^{\bullet}(M_1,\ldots,M_n)$
as the total complex associated with a bicomplex, where the
bidegree $(\deg_0,\deg_1)$ on $M_1\otimes T(R_+)\otimes\ldots\otimes
T(R_+)\otimes M_n$ is given by
$$\deg_0(x_1\otimes t_1\otimes\ldots\otimes t_{n-1}\otimes x_n)=
\sum_{i\equiv 0 (2)}\deg(t_i),$$
$$\deg_1(x_1\otimes t_1\otimes\ldots\otimes t_{n-1}\otimes x_n)=
\sum_{i\equiv 1 (2)}\deg(t_i),$$
where $t_i\in T(R_+)$, $x_i\in M_i$, $\deg$ denotes the cohomological
degree on $T(R_+)$.
Therefore, there is a spectral sequence abbuting to cohomology of
$B^{\bullet}(M_1,\ldots,M_n)$ with the $E_1$-term
$$H^*(M_1\otimes T(R_+)\otimes M_2)\otimes T(R_+)\otimes
H^*(M_3\otimes T(R_+)\otimes M_4)\otimes\ldots,$$
where the last factor of the tensor product
is either $M_n$ or $H^*(M_{n-1}\otimes T(R_+)\otimes M_n)$.
Using part (ii) we see that $E_1$ is isomorphic to the complex
of the form
$$B^{\bullet}(M_1',\ldots,M'_{n'})[(n-n')(d+1)]$$
with $n'<n$. It remains to apply the induction assumption.

\noindent
(iv) Consider first the case $n=1$. The complex
$B^{\bullet}(k,M_1,k)=T(R_+)\otimes M_1\otimes T(R_+)$ is the total
complex of the bicomplex $(\partial_1\otimes\id, \id\otimes\partial_2)$,
where $\partial_1$ and $\partial_2$ are bar-differentials on
$T(R_+)\otimes M_1$ and $M_1\otimes T(R_+)$. Our assertion follows
immediately from (i) by considering the spectral sequence
associated with this bicomplex.

Now assume that for some $n>1$
the assertion holds for all $n'<n$. As before we view
$B^{\bullet}(k,M_1,\ldots,M_n,k)$ as the total complex of a bicomplex
by defining the bidegree 
on $T(R_+)\otimes M_1\otimes\ldots\otimes M_n
\otimes T(R_+)$ as follows:
$$\deg_0(t_0\otimes x_1\otimes t_1\ldots\otimes x_n\otimes t_n)=
\sum_{i\equiv 0 (2)}\deg(t_i),$$
$$\deg_1(t_0\otimes x_1\otimes t_1\ldots\otimes x_n\otimes t_n)=
\sum_{i\equiv 1 (2)}\deg(t_i).$$
Assume first that $n$ is even.
Then there is a spectral sequence associated with the above bicomplex
abutting to
the cohomology of $B^{\bullet}(k,M_1,\ldots,M_n,k)$ and with the $E_1$-term
$$T(R_+)\otimes H^*(M_1\otimes T(R_+)\otimes M_2)\otimes T(R_+)
\otimes\ldots \otimes H^*(M_{n-1}\otimes T(R_+)\otimes M_n)\otimes T(R_+).$$
Using (ii) we see that $E_1$ is isomorphic
to the complex of the form
$$B^{\bullet}(k,M'_1,\ldots,M'_{n/2},k)[n(d+1)/2],$$
so we can apply the induction assumption.
If $n$ is odd then we consider another spectral sequence associated
with the above bicomplex, so that
\begin{align*}
E_1=&H^*(T(R_+)\otimes M_1)\otimes T(R_+)\otimes\\
&H^*(M_2\otimes T(R_+)\otimes M_3)\otimes T(R_+)
\otimes\ldots \otimes H^*(M_{n-1}\otimes T(R_+)\otimes M_n)\otimes T(R_+).
\end{align*}
By (i) and (ii) this complex is isomorphic to
$B^{\bullet}(k,M'_1,\ldots,M'_{(n-1)/2},k)[(n+1)(d+1)/2]$.
Again we can finish the proof by applying the induction assumption.
\ed


We will also need the following simple lemma.

\begin{lem}\label{convlem}
Let $C^{\bullet}$ be a complex in an abelian category
equipped with a decreasing filtration $C^{\bullet}=F^0 C^{\bullet}\supset
F^1 C^{\bullet}\supset F^2 C^{\bullet}\supset\ldots$
such that $C^n=\limproj_i C^n/F^i C^n$ for all $n$.
Let $\gr_i C^{\bullet}= F^i C^{\bullet}/F^{i+1} C^{\bullet}$, $i=0,1,\ldots$
be the associated graded factors.
Assume that $H^n \gr_i C^{\bullet}=0$ for all $i>0$ and for some fixed $n$.
Then the natural map $H^n C^{\bullet}\ra H^n \gr_0 C^{\bullet}$ is
injective.
\end{lem}

\Pf . Considering an exact sequence of complexes
$$0\ra F^1 C^{\bullet}\ra C^{\bullet}\ra \gr_0 C^{\bullet}\ra 0$$
one can easily reduce the proof to the case
$H^n \gr_i C^{\bullet}=0$ for all $i\ge 0$. In this case we have to show
that $H^n C^{\bullet}=0$. Let $c\in C^n$ be a cocycle and let
$c_i$ be its image in $C^n/F^i C^n$. It suffices to construct
a sequence of elements $x_i\in C^{n-1}/F^i C^{n-1}$, $i=1,2,\ldots$, such that
$x_{i+1}\equiv x_i\mod F^i C^{n-1}$ and $c_i=d(i)x_i$, where $d(i)$
is the differential on $C^{\bullet}/F^i C^{\bullet}$.
Since $n$-th cohomology of $C^{\bullet}/F^1 C^{\bullet}=\gr_0 C^{\bullet}$
is trivial we can find $x_1$ such that $c_1=d(1)x_1$. Then we proceed
by induction: once $x_1,\ldots,x_i$ are chosen
an easy diagram chase using the exact triple of complexes
$$0\ra \gr_i C^{\bullet}\ra C^{\bullet}/F^{i+1} C^{\bullet}\ra
C^{\bullet}/F^i C^{\bullet}\ra 0$$
and the vanishing of $H^n(\gr_i C^{\bullet})$ show that $x_{i+1}$ exists.
\ed

\begin{thm}\label{vanishthm}
Under the assumptions of Theorem \ref{mainthm}
one has $HH^i_{0,md}(A)=0$ for $i<m(d+2)$ and
$\dim HH^{d+2}_{0,d}(A)\le 1$, where $A=A_L$.
\end{thm}

\Pf . Set $C^i=C^i_{0,md}(A)$ (see \ref{obstrsec}). 
Note that Hochschild differential
maps $C^i$ into $C^{i+1}$ (since $m_2$ preserves both gradings on $A$).
Recall that the decomposition of $A$ into
graded pieces with respect to the cohomological degree has form
$A=R\oplus M$, where $R$ has degree $0$ and
$M=\oplus_{i\in\Z}H^d(X,L^i)$ has degree $d$.
The natural augmentation of $A$ is given by the ideal $A_+=R_+\oplus M$.
Each of the spaces $C^i$ decomposes into a direct sum
$C^i=C^i(0)\oplus C^i(d)$, where 
$C^i(0)\sub\Hom(A_+^{\otimes i},R)$,
$C^i(d)\sub\Hom(A^{\otimes i}_+,M)$.
More precisely, the space $C^i(0)$
consists of linear maps 
\begin{equation}\label{Ci0eq}
[T(R_+)\otimes M\otimes T(R_+)\otimes\ldots\otimes M\otimes T(R_+)]_i\ra R
\end{equation}
preserving the internal grading,
where there are $m$ factors of $M$ in the tensor product and the index
$i$ refers to the total number of factors $H^*(L^*)$ (so that the LHS can be
considered as a subspace of $A_+^{\otimes i}$).
Similarly, the space $C^i(d)$ consists of linear maps
$$[T(R^+)\otimes M\otimes T(R_+)\otimes\ldots\otimes M\otimes T(R_+)]_i
\ra M$$
preserving the internal grading,
where there is $m+1$ factors of $M$ in the tensor product. 
Clearly, $C^{\bullet}(d)$ is a subcomplex in $C^{\bullet}$, so we
have an exact sequence of complexes
$$0\ra C^{\bullet}(d)\ra C^{\bullet}\ra C^{\bullet}(0)\ra 0.$$
Therefore, it suffices to prove that 
$H^i(C^{\bullet}(0))=H^i(C^{\bullet}(d))=0$ for $i<m(d+2)$,
and that in the case $m=1$ one has in addition
$H^{d+2}(C^{\bullet}(d))=0$ and $\dim H^{d+2}(C^{\bullet}(0))\le 1$.

To compute the cohomology of these two complexes we can use
spectral sequences associated with some natural filtrations to 
reduce the problem to simpler complexes.
First, let us consider the decomposition
$$C^{\bullet}(0)=\prod_{j\ge 0}C^{\bullet}(0)_j,$$
where $C^i(0)_j\sub C^i(0)$ is the space of maps (\ref{Ci0eq})  
with the image contained in $H^0(L^j)\sub R$.
The differential on $C^{\bullet}(0)$ has form
$$\de(x_j)_{j\ge 0}=(\sum_{j'\le j}\de_{j',j} x_{j'})_{j\ge 0}$$
for some maps $\de_{j',j}:C^{\bullet}(0)_{j'}\ra C^{\bullet}(0)_j$,
where $j'\le j$. By Lemma \ref{convlem} it suffices to prove
that one has
$H^i(C^{\bullet}(0)_j,\de_{j,j})=0$ for $i<m(d+2)$ and all $j$,
while for $m=1$ one has in addition
$H^{d+2}(C^{\bullet}(0)_j,\de_{j,j})=0$ for $j>0$
and $\dim H^{d+2}(C^{\bullet}(0)_0,\de_{0,0})\le 1$.
But
$$(C^{\bullet}(0)_j,\de_{j,j})=\Hom(K^{\bullet}_{m,j},R_j)[-m],$$
where $K_m^{\bullet}=B^{\bullet}(k,M,\ldots,M,k)$
($m$ copies of $M$) and $K^{\bullet}_{m,j}$ is its $j$-th
graded component with respect to the internal
grading. Here we use the following convention for the grading on the
dual complex: $\Hom(K^{\bullet},R)^i=\Hom(K^{-i},R)$.
Therefore, Proposition \ref{cohprop}(iv) implies
that cohomology of $C^{\bullet}(0)_j$
is concentrated in degrees $\ge m(d+1)+m=m(d+2)$. 
Moreover, for $m=1$ the $(d+2)$-th cohomology space is non-zero
only for $j=0$, in which case it is one-dimensional.

For the complex $C^{\bullet}(d)$ we have to use a different
filtration (since $M$ is not bounded below with respect to
the internal grading). Consider the decreasing filtration on $C^{\bullet}(d)$
induced by the following grading on 
$T(R^+)\otimes M\otimes T(R_+)\otimes\ldots\otimes M\otimes T(R_+)$:
$$\deg(t_1\otimes x_1\otimes t_2\otimes\ldots\otimes x_{m+1}\otimes t_{m+2})=
\deg(t_1)+\deg(t_{m+2}),$$
where $t_i\in T(R_+)$, $x_i\in M$, the degree of $R_+$ is taken to be $-1$.
The associated graded complex is
$$\Hom_{\gr}(T(R^+)\otimes B^{\bullet}(M,\ldots,M)\otimes T(R^+), M)[-m-1],$$
where there are $m+1$ factors of $M$ in the bar-construction.
It remains to apply Proposition \ref{cohprop}(iii).
\ed


\subsection{Some Massey products}

In this subsection we will show the nontriviality of the
canonical class of $A_{\infty}$-structures 
on $A_L$ and combine it with our computations of the Hochschild
cohomology to prove the main theorem.

Note that the canonical class of $A_{\infty}$-structures can
be defined in a more general context. Namely, if $\CC$ is an abelian
category with enough injectives then we can define the canonical class
of $A_{\infty}$-structures on the derived category $D^+(\CC)$ of bounded
below complexes over $\CC$. Indeed, one can use the equivalence of
$D^+(\CC)$ with the homotopy category of complexes with injective terms
and apply Kadeishvili's construction to the dg-category of such complexes
(see \cite{P-b},1.2 for more details). In the case when $\CC$
is the category of coherent sheaves the resulting $A_{\infty}$-structure
is strictly $A_{\infty}$-isomorphic to the structure obtained using
{\cech} resolutions (since the relevant dg-categories are equivalent). 
In this context we have the following construction
of nontrivial Massey products. 

\begin{lem}\label{MPlem} Let
$\CC$ be an abelian category with enough injectives,
$$0\ra \FF_0\stackrel{\a_1}{\ra}\FF_1
\stackrel{\a_2}{\ra}\FF_2\ra\ldots\stackrel{\a_n}{\ra}\FF_n \ra 0$$
be an exact sequence in $\CC$, where $n\ge 2$,
and let $\b:\FF_n\ra\FF_0[n-1]$
be a morphism in the derived category $D^b(\CC)$ corresponding to
the Yoneda extension class in $\Ext^{n-1}(\FF_n,\FF_0)$ represented
by the above sequence. 
Assume that $\Ext^{j-i-1}(\FF_j,\FF_i)=0$ when $0\le i<j\le n-1$. Then
$$m_{n+1}(\a_1,\ldots,\a_n,\b)=\pm\id_{\FF_0}$$
for any $A_{\infty}$-structure $(m_i)$ on $D^b(\CC)$ from
the canonical class.
\end{lem}

\Pf . Assume first that $n=2$. Then $m_3(\a_1,\a_2,\b)$ is the unique value
of the well-defined triple Massey product in $D^b(\CC)$ (see
\cite{P-AYBE}, 1.1).
Using the standard recipe for its calculation (see \cite{GM}, IV.2)
we immediately get that $m_3(\a_1,\a_2,\b)=\id$.

For general $n$ we can proceed by induction. Assume that the statement 
is true for $n-1$. Set $\FF'_{n-1}=\ker(\a_n)$. Then we have
exact sequences
\begin{equation}\label{seq1}
0\ra\FF_0\stackrel{\a_1}{\ra}\FF_1\ra\ldots\ra
\FF_{n-2}\stackrel{\a'_{n-1}}{\ra}\FF'_{n-1}\ra 0,
\end{equation}
\begin{equation}\label{seq2}
0\ra\FF'_{n-1}\stackrel{i}{\ra}\FF_{n-1}\stackrel{\a_n}{\ra}\FF_n\ra 0,
\end{equation}
where $m_2(\a'_{n-1},i)=i\circ\a'_{n-1}=\a_{n-1}$. By the definition, we have
$\b=m_2(\ga,\b')$, where $\b'\in\Ext^{n-2}(\FF'_{n-1},\FF_0)$ and
$\ga\in\Ext^1(\FF_n,\FF'_{n-1})$ are the extension classes corresponding
to these exact sequences. Applying the $A_{\infty}$-axiom
to the elements $(\a_1,\ldots,\a_n,\ga,\b')$ and using the
vanishing of $m_{n-i+2}(\a_{i+1},\ldots,a_n,\ga,\b')\in
\Ext^{i-1}(\FF_i,\FF_0)$ for $0<i<n$, we get
\begin{equation}\label{maseq1}
m_{n+1}(\a_1,\ldots,\a_n,\b)=
\pm m_n(\a_1,\ldots,\a_{n-1},m_3(\a_{n-1},\a_n,\ga),\b').
\end{equation}
Furthermore, applying the $A_{\infty}$-axiom to
the elements $(\a'_{n-1},i,\a_n,\ga)$ we get
\begin{equation}\label{maseq2}
m_3(\a_{n-1},\a_n,\ga)=\pm m_2(\a'_{n-1},m_3(i,\a_n,\ga)).
\end{equation}
Next, we claim that sequences (\ref{seq1}) and (\ref{seq2}) satisfy
the assumptions of the lemma. Indeed, for (\ref{seq1})
this is clear, so we just have to check that
$\Hom(\FF_{n-1},\FF'_{n-1})=0$. The exact sequence (\ref{seq1})
gives a resolution $\FF_0\ra\ldots\ra\FF_{n-2}$ of $\FF'_{n-1}$ and
we can compute $\Hom(\FF_{n-1},\FF'_{n-1})$ using this resolution.
Now the required vanishing follows from our assumption that
$\Ext^{n-i-2}(\FF_{n-1},\FF_i)=0$ for $0\le i\le n-2$.
Therefore, we have
$$m_3(i,\a_n,\ga)=\id.$$
Together with (\ref{maseq2}) this implies that
$$m_3(\a_{n-1},\a_n,\ga)=\pm\a'_{n-1}.$$
Substituting this into (\ref{maseq1}) we get
$$m_{n+1}(\a_1,\ldots,\a_n,\b)=
\pm m_n(\a_1,\ldots,\a_{n-1},\a'_{n-1},\b').$$
It remains to apply the induction assumption to the sequence
(\ref{seq1}).
\ed

\noindent
{\it Proof of Theorem \ref{mainthm}.}
(i) Since the algebra $A=A_L$ is concentrated in degrees $0$ and $d$,
the first potentially nontrivial higher product of an admissible
$A_{\infty}$-structure $(m_i)$ on $A$ is $m_{d+2}$.
Therefore, by Lemma \ref{obstrlem} for every such $A_{\infty}$-structure
$(m_i)$ on $A$ the map $m_{d+2}$ induces a cohomology class
$[m_{d+2}]\in HH^{d+2}_{0,d}(A)$. We claim that if $(m'_i)$ is
another admissible $A_{\infty}$-structure on $A$ then
$(m_i)$ is strictly $A_{\infty}$-isomorphic 
to $(m'_i)$ if and only if $[m_{d+2}]=[m'_{d+2}]$.
Indeed, this follows from Lemma \ref{obstrlem} and from the
vanishing of higher obstructions due to Theorem \ref{vanishthm}
(these obstructions lie in $HH^{md+2}_{0,md}(A)$ where $m>1$, and
the vanishing follows since $md+2<m(d+2)$).
In particular, an admissible $A_{\infty}$-structure $(m_i)$ is nontrivial
if and only if $[m_{d+2}]\neq 0$.
Since by Theorem \ref{vanishthm}
the space $HH^{d+2}_{0,d}(A)$ is at most one-dimensional, it
remains to prove the nontriviality of an admissible $A_{\infty}$-structure
from the canonical class. Replacing $L$ by its sufficiently high
power if necessary we can assume that there exists $d+1$ sections
$s_1,\ldots,s_{d+1}\in H^0(L)$ without common zeroes. The corresponding
Koszul complex gives an exact sequence
$$0\ra\OO\ra \OO^{\oplus (d+1)}\otimes_{\OO} L\ra
\OO^{\oplus {d+1\choose 2}}\otimes_{\OO} L^2\ra\ldots\ra
\OO^{\oplus (d+1)}\otimes_{\OO} L^d\ra L^{d+1}\ra 0.$$
By our assumptions this sequence satisfies the conditions required
in Lemma \ref{MPlem}, hence we get a nontrivial
$(d+2)$-ple Massey product for our $A_{\infty}$-structure.

\noindent (ii) Applying Lemma \ref{obstrlem2} we see that obstructions
for connecting two strict $A_{\infty}$-isomorphisms by a homotopy lie
in $\oplus_{m\ge 1}HH^{md+1}_{0,md}(A)$. But this space is zero by
Theorem \ref{vanishthm}.
\ed

\begin{cor} Under assumptions of Theorem \ref{mainthm}
the space $HH^{d+2}_{0,d}(A_L)$ is one-dimensional.
\end{cor}

\Pf . Indeed, from Theorem \ref{vanishthm} we know that 
$\dim HH^{d+2}_{0,d}(A_L)\le 1$.
If this space were zero then the above argument would
show that all admissible $A_{\infty}$-structures on $A_L$ are trivial.
But we know that $A_{\infty}$-structures on $A_L$ from the
canonical class are nontrivial.
\ed

\begin{rem} One can ask
whether there exists an $A_{\infty}$-structure on $A_L$ from the
canonical class such that $m_n=0$ for $n>d+2$ or at least $m_n=0$ for
all sufficiently large $n$. However, even in
the case of smooth curves of genus $\ge 1$
the answer is ``no''. The proof can be obtained using the
construction of a universal deformation of a coherent sheaf (when it exists)
using the canonical $A_{\infty}$-structure, outlined in \cite{P-b}. For example,
it is shown there that the products 
$$m_{n+2}:H^1(\OO_X)^{\otimes n}\otimes H^0(L^{n_1})\otimes H^0(L^{n_2})\ra
H^0(L^{n_1+n_2})$$
appear as coefficients in the universal deformation of the structure sheaf.
The base of this family is $\Spec R$, where $R\simeq k[[t_1,\ldots,t_g]]$ 
is the completed symmetric algebra of $H^1(\OO_X)^{\vee}$. If all sufficiently
large products were zero, this family would be induced by the base change from some family over
an open neighborhood $U$ of zero in the affine space $\A^g$. 
But this would imply that the
embedding of $\Spec R$ into the Jacobian (corresponding to the isomorphism of
$R$ with the completion of the local ring of the Jacobian at zero) 
factors through $U$, which is false.
\end{rem}

\subsection{Proof of Theorem \ref{groupthm}}

Theorem \ref{mainthm}(i) easily implies that every admissible
$A_{\infty}$-structure on $A=A_L$ is (strictly)
$A_{\infty}$-isomorphic to some strictly
unital $A_{\infty}$-structure. Therefore, it
is enough to prove our statement for strictly unital structures.
Recall that the group of strict $A_{\infty}$-isomorphisms $HG$ is the group of
coalgebra automorphisms of $\Bar(A_L)$ inducing the identity
map $A_L\ra A_L$ and preserving two grading on $\Bar(A_L)$
induced by the two gradings of $A_L$.
Thus, we can identify $HG$ with a subgroup
of algebra automorphisms of the completed cobar-construction
$\Cobar(A_L)=\prod_{n\ge 0}T^n(A_L^*[-1])$
(our convention is that passing to dual vector space changes the grading to
the opposite one).

By Theorem \ref{mainthm}(ii) for every strict 
$A_{\infty}$-automorphism $f$ of an
$A_{\infty}$-structure $m$ there exists a homotopy
from $f$ to the trivial $A_{\infty}$-automorphism $f^{tr}$. Let
$\a=\a_f^*$ be the automorphism of $\Cobar(A_L)$ corresponding to $f$
and $h=H^*:\Cobar(A_L)\ra\Cobar(A_L)[-1]$ be the map giving the
homotopy from $f$ to $f^{tr}$.
The equations dual to (\ref{hommor1}) and (\ref{hommor2}) in our
case have form
$$h(xy)=h(x)y\pm\a(x)h(y),$$
$$\a=\id+d\circ h+h\circ d,$$
where $d$ is the differential on $\Cobar(A_L)$ associated with $m$.
Recall that $A_L=H^0\oplus H^1$, where $H^0=\oplus_{n\ge 0} H^0(X,L^n)$,
$H^1=\oplus_{n\le 0} H^1(X,L^n)$.
Since $h$ has degree $-1$ 
we have $h((H^1)^*[-1])=0$ and $h((H^0)^*[-1])\sub\hat{T}((H^1)^*[-1])$.
Furthermore, since $h$ preserves the internal degree, we have
$h(H^0(X,L^n)^*[-1])=0$ for all $n>0$. 
Let $\eps\in (H^0)^*[-1]\sub\Cobar(A_L)$ be an element corresponding
to the natural
projection $H^0\ra H^0(X,\OO_X)\simeq k$. Then
we have $A_L^*[-1]=k\eps\oplus V$, where
$V=(H^1)^*[-1]\oplus(\oplus_{n>0}H^0(X,L^n))^*[-1]$, and
$h(V)=0$. Let $\lan V\ran\sub\Cobar(A_L)$ be the subalgebra
topologically generated by $V$. Then $h$ vanishes on $\lan V\ran$.
Also, for every $x\in V$ we have
$$dx=\eps x+ x\eps\mod\lan V\ran$$
since our $A_{\infty}$-structure is strictly unital.
Hence, for $x\in V$ we have
$$\a(x)=x+h(dx)=x+h(\eps)x-\a(x)h(\eps),$$
which implies that
$$\a(x)=(1+h(\eps))x(1+h(\eps))^{-1}.$$
In particular, the restriction of $\a$ to the subalgebra
$\hat{T}(H^1(X,\OO_X)^*[-1])$ is the inner automorphism
associated with the invertible
element $1+h(\eps)\in\hat{T}(H^1(X,\OO_X)^*[-1])$.
On the other hand, we have
$$d\eps=\eps^2\mod\lan V\ran.$$
Hence,
$$\a(\eps)=\eps+dh(\eps)+h(\eps)\eps-\a(\eps)h(\eps),$$
so that
$$\a(\eps)=(1+h(\eps))\eps(1+h(\eps))^{-1}+dh(\eps)\cdot(1+h(\eps))^{-1}.$$
Thus, $\a$ is uniquely determined by $h(\eps)$. Also, by Lemma 
\ref{closmorlem} $h(\eps)$ can
be an arbitrary element of 
$\prod_{n\ge 1}T^n(H^1(X,\OO_X)^*[-1])$.
\ed

\end{document}